\documentclass[12 pt]{article}
\hoffset - 10 mm

\catcode`\@=11 \@addtoreset{equation}{section}

\catcode`\@=12
\usepackage{color}
\newcommand{\mathbb}{\mathbf}

\newtheorem{Theorem}{Theorem}[section]
\newtheorem{Proposition}{Proposition}[section]
\newtheorem{Lemma}{Lemma}[section]
\newtheorem{Corollary}{Corollary}[section]

\newcommand{\bTheorem}[1]{
\begin{Theorem} \label{T#1} }
\newcommand{\eT}{\end{Theorem}}

\newcommand{\bProposition}[1]{
\begin{Proposition} \label{P#1}}
\newcommand{\eP}{\end{Proposition}}

\newcommand{\bLemma}[1]{
\begin{Lemma} \label{L#1} }
\newcommand{\eL}{\end{Lemma}}

\newcommand{\bCorollary}[1]{
\begin{Corollary} \label{C#1} }
\newcommand{\eC}{\end{Corollary}}

\newcommand{\bFormula}[1]{
\begin{equation} \label{#1}}
\newcommand{\eF}{\end{equation}}

\newcommand{\Ov}[1]{\overline{#1}}
\newcommand{\DC}{C^\infty_c}

\newcommand{\vr}{\varrho}

\newcommand{\vt}{\vartheta}
\newcommand{\vu}{\vc{u}}
\newcommand{\vc}[1]{{\bf #1}}

\newcommand{\Div}{{\rm div}_x}
\newcommand{\Grad}{\nabla_x}

\newcommand{\tn}[1]{\mbox {\F #1}}
\newcommand{\dx}{{\rm d} x}
\newcommand{\tvV}{{\tilde \vc{V}}}
\newcommand{\dt}{{\rm d} t }

\newcommand{\tvr}{\tilde \vr}
\newcommand{\tvt}{\tilde \vt}
\newcommand{\tvu}{\tilde \vu}

\newcommand{\intO}[1]{\int_{\Omega} #1 \ \dx}

\newcommand{\ep}{\varepsilon}

\font\F=msbm10 scaled 1000

\definecolor{Cgrey}{rgb}{0.85,0.85,0.85}
\definecolor{Cblue}{rgb}{0.50,0.85,0.85}
\definecolor{Cred}{rgb}{0.80,0.50,0.50}
\definecolor{fancy}{rgb}{0.10,0.85,0.10}

\newcommand\Cbox[2]{%
    \newbox\contentbox%
    \newbox\bkgdbox%
    \setbox\contentbox\hbox to \hsize{%
        \vtop{
            \kern\columnsep
            \hbox to \hsize{%
                \kern\columnsep%
                \advance\hsize by -2\columnsep%
                \setlength{\textwidth}{\hsize}%
                \vbox{
                    \parskip=\baselineskip
                    \parindent=0bp
                    #2
                }%
                \kern\columnsep%
            }%
            \kern\columnsep%
        }%
    }%
    \setbox\bkgdbox\vbox{
        \color{#1}
        \hrule width  \wd\contentbox %
               height \ht\contentbox %
               depth  \dp\contentbox
        \color{black}
    }%
    \wd\bkgdbox=0bp%
    \vbox{\hbox to \hsize{\box\bkgdbox\box\contentbox}}%
    \vskip\baselineskip%
}

\date{}
\begin{document}


\title{A regularity criterion for the weak solutions to the Navier-Stokes-Fourier system}
\author{Eduard Feireisl \thanks{Eduard Feireisl acknowledges the support of the project LL1202 in the
programme ERC-CZ funded
by the Ministry of Education, Youth and Sports of the Czech Republic.} \and Anton\' \i n Novotn\' y \and Yongzhong Sun \thanks{Yongzhong Sun is supported by NSF of China(Grants No. 11171145 and 10931007) and the PAPD of
Jiangsu Higher Education Institutions.}}

\maketitle

\bigskip

\centerline{Charles University in Prague, Faculty of Mathematics and Physics, Mathematical Institute}
\centerline{Sokolovsk\' a 83, 186 75 Prague 8,
Czech Republic}

\medskip

\centerline{IMATH Universit\' e du Sud Toulon-Var}
\centerline{BP 132, 83957 La Garde, France}

\medskip

Department of Mathematics, Nanjing University,
Nanjing, Jiangsu 210093, China

\medskip

\begin{abstract}

We show that any weak solution to the full Navier-Stokes-Fourier system emanating from the data belonging to the Sobolev space
$W^{3,2}$ remains regular as long as the velocity gradient is bounded. The proof is based on the weak-strong uniqueness property
and parabolic \emph{a priori} estimates for the local strong solutions.

\end{abstract}

{\bf Key words:} Navier-Stokes-Fourier system, weak solution, regularity

\bigskip


\bigskip

\tableofcontents

\vglue 2 cm

\section{Introduction}
\label{i}

In continuum mechanics, the motion of a general compressible, viscous, and heat conducting fluid is described by the thermostatic variables -
the mass density $\vr = \vr(t,x)$ and the absolute temperature $\vt(t,x)$, and the velocity field $\vu = \vu(t,x)$ evaluated at the time
$t$ and the spatial position $x$ belonging to the reference physical domain $\Omega \subset R^3$. The time evolution of the fluid, emanating from the initial data
\bFormula{i1}
\vr(0,\cdot) = \vr_0, \ \vt(0, \cdot) = \vt_0,\ \vu(0, \cdot) = \vu_0 \ \mbox{in} \ \Omega,
\eF
is governed by the \emph{Navier-Stokes-Fourier} system of partial differential equations:

\Cbox{Cgrey}{

\centerline{\textsc{Equation of continuity:}}
\bFormula{i2}
\partial_t \vr + \Div (\vr \vu) = 0;
\eF

\centerline{\textsc{Momentum equation:}}
\bFormula{i3}
\partial_t (\vr \vu) + \Div (\vr \vu \otimes \vu) + \Grad p(\vr, \vt) = \Div \tn{S}(\vt, \Grad \vu);
\eF

\centerline{\textsc{Total energy balance:}}
\bFormula{i4}
\partial_t \left( \frac{1}{2} \vr |\vu|^2 + \vr e(\vr, \vt) \right) +
\Div \left[ \left( \frac{1}{2} \vr |\vu|^2 + \vr e(\vr, \vt) \right) \vu \right] + \Div (p(\vr, \vt) \vu)
\eF
\[
= \Div \left( \tn{S}(\vt, \Grad \vu) \vu \right) - \Div \vc{q}(\vt, \Grad \vt).
\]

}

The symbols $p = p(\vr,\vt)$ and $e=e(\vr, \vt)$ stand for the \emph{pressure} and the (specific) \emph{internal energy}, respectively. Furthermore,
$\tn{S} = \tn{S}(\vt, \Grad \vu)$ denotes the \emph{viscous stress} given by

\Cbox{Cgrey}{

\centerline{\textsc{Newton's rheological law:}}
\bFormula{i5}
\tn{S}(\vt, \Grad \vu) = \mu(\vt) \left( \Grad \vu + \Grad^t \vu - \frac{2}{3} \Div \vu \tn{I} \right),
\eF

}

\noindent while $\vc{q} = \vc{q}(\vt, \Grad \vt)$ is the \emph{heat flux} determined by

\Cbox{Cgrey}{

\centerline{\textsc{Fourier's law:}}
\bFormula{i6}
\vc{q}(\vt, \Grad \vu) = - \kappa(\vt) \Grad \vt.
\eF

}

We restrict ourselves to \emph{bounded} regular domains $\Omega \subset R^3$, with energetically insulated boundary. More specifically, we
consider the standard

\Cbox{Cgrey}{

\centerline{\textsc{No slip boundary condition:}}
\bFormula{i7}
\vu|_{\partial \Omega} = 0,
\eF

}

\noindent together with

\Cbox{Cgrey}{

\centerline{\textsc{No-flux condition:}}
\bFormula{i8}
\vc{q} \cdot \vc{n}|_{\partial \Omega} = 0,
\eF

}

\noindent where $\vc{n}$ is the (outer) normal vector to $\partial \Omega$. Under these circumstance, the total mass as well as the total energy of the system are constants of motion, specifically, integrating equations (\ref{i2}), (\ref{i4}) over $\Omega$ we obtain
\bFormula{i9}
\intO{ \vr(t, \cdot) } = \intO{ \vr_0 } = M_0,
\eF
\bFormula{i10}
\intO{ \left( \frac{1}{2} \vr |\vu|^2 + \vr e(\vr, \vt) \right) (t, \cdot) } =
\intO{ \left( \frac{1}{2} \vr_0 |\vu_0|^2 + \vr_0 e(\vr_0, \vt_0) \right)  } = E_0.
\eF

From the mathematical viewpoint, the Navier-Stokes-Fourier system suffers the same deficiency as most of its counterparts in continuum mechanics -
the lack of sufficiently strong {\it a priori} bounds. It is known that the problem (\ref{i1} - \ref{i8}) is well posed, locally in time, in the
framework of classical solutions, see Tani \cite{TAN}, or in slightly more general \emph{energy spaces} of Sobolev type, see
Valli \cite{Vall2}, \cite{Vall1}, Valli and Zajaczkowski,
\cite{VAZA}, among others. Besides these local existence results, Matsumura and Nishida established in a series of papers \cite{MATS}, \cite{MANI1}, \cite{MANI} global-in-time existence of strong/classical solutions provided the initial data are sufficiently close to an equilibrium solution.
Similarly to the well known \emph{incompressible} Navier-Stokes system, the existence of global-in-time smooth solutions for any physically admissible and large initial data remains an outstanding open problem.

The system (\ref{i1} - \ref{i8}) has been also studied in the framework of \emph{weak solutions}. Hoff and Jenssen \cite{HOJE} established global existence for radially symmetric data in $R^3$. They also identified one of the main stumbling blocks in the analysis of the
Navier-Stokes-Fourier system, namely the (hypothetical) appearance of \emph{vacuum zones}, where the density vanishes and the classical understanding of the equations breaks down. More recently, Bresch and Desjardins \cite{BRDE}, \cite{BRDE1} discovered a new {\it a priori} bound on the density gradient leading to global-in-time existence in the truly $3D-$setting conditioned, unfortunately, by a very specific relation satisfied by
the density dependent viscosity coefficients and a rather unrealistic formula for the pressure that
has to be infinite (negative) for $\vr \to 0$.

\subsection{Weak and dissipative solutions}

We adopt the approach to the Navier-Stokes-Fourier system originated in \cite{EF73} and further developed and detailed in \cite{FeNo6}.
Suppose that $p$ and $e$ are interrelated through

\centerline{\textsc{Gibbs' equation:}}
\bFormula{i11}
\vt Ds(\vr, \vt) = D e (\vr, \vt) + p(\vr, \vt) D \left( \frac{1}{\vr} \right),
\eF
where $s = s(\vr, \vt)$ is a new thermodynamic function called (specific) \emph{entropy}. If $\vr$, $\vt$ are smooth and bounded below away from zero and if $u$ is smooth, then the total energy balance (\ref{i4}) can be replaced by

\centerline{\textsc{Thermal energy balance:}}
\bFormula{i12}
\vr \frac{\partial e(\vr, \vt)}{\partial \vt} \left( \partial_t \vt + \vu \cdot \Grad \vt \right) - \Div (\kappa(\vt) \Grad \vt) =
\tn{S}(\vt, \Grad \vu) : \Grad \vu - \vt \frac{\partial p(\vr, \vt)}{\partial \vt} \Div \vu.
\eF

Furthermore, dividing (\ref{i12}) on $\vt$ we arrive at

\centerline{\textsc{Entropy (production) equation:}}
\bFormula{i13}
\partial_t (\vr s(\vr, \vt)) + \Div (\vr s(\vr, \vt) \vu) + \Div \left( \frac{\vc{q}(\vt, \Grad \vt)}{\vt} \right) = \sigma,
\eF
with

\centerline{\textsc{Entropy production rate:}}
\bFormula{i14}
\sigma = \frac{1}{\vt} \left( \tn{S}(\vt, \Grad \vu) : \Grad \vu - \frac{\vc{q} \cdot \Grad \vt}{\vr} \right).
\eF
Let us remark that the systems (\ref{i2}, \ref{i3}, \ref{i4}) , (\ref{i2}, \ref{i3}, \ref{i12}), and (\ref{i2}, \ref{i3}, \ref{i14})
are perfectly \emph{equivalent} in the class of classical solutions.

Another crucial observation, which is the platform of the approach developed in \cite{FeNo6}, is that we can relax (\ref{i14}) to
\bFormula{i15}
\sigma \geq \frac{1}{\vt} \left( \tn{S}(\vt, \Grad \vu) : \Grad \vu - \frac{\vc{q} \cdot \Grad \vt}{\vr} \right)
\eF
provided the system is supplemented by the integrated total energy balance (\ref{i10}), meaning the relations (\ref{i2}, \ref{i3}, \ref{i10},
\ref{i13}), with (\ref{i15}) are still equivalent to the original system  (\ref{i2}, \ref{i3}, \ref{i4}), at least for smooth solutions.
The resulting problem is mathematically tractable, more specifically, we report the following results:

\begin{itemize}
\item {\textsc{Existence.}} The problem \{ (\ref{i1} - \ref{i3}), (\ref{i5} - \ref{i10}), (\ref{i13}), (\ref{i15}) \} admits global-in-time
weak solutions for any finite energy initial data under certain structural restrictions imposed on $p$, $e$, $s$, and $\mu$, $\kappa$
presented in Section \ref{m} below, see \cite[Theorem 3.1]{FeNo6}.
\item {\textsc{Compatibility.}} Any weak solution of \{ (\ref{i1} - \ref{i3}), (\ref{i5} - \ref{i10}), (\ref{i13}), (\ref{i15}) \} that is
regular solves the original problem (\ref{i1} - \ref{i8}), see \cite[Chapter 2]{FeNo6}.
\item {\textsc{Weak-strong uniqueness.}} Any weak solution coincides with a classical solution emanating from the same initial data provided the latter exists, see \cite{FeiNov10}.

\end{itemize}

\subsubsection{Relative entropy and dissipative solutions}

The property of \emph{weak-strong uniqueness} mentioned above is closely related to the \emph{relative entropy inequality} introduced in
\cite{FeiNov10}. We start by introducing another thermodynamic function:

\centerline{\textsc{Ballistic free energy:}}
\bFormula{i16}
H_\Theta(\vr, \vt) = \vr e(\vr, \vt) - \Theta \vr s(\vr, \vt),
\eF
see \cite{Eri}, together with

\centerline{\textsc{Relative entropy functional:}}
\bFormula{i17}
\mathcal{E} \left( \vr, \vt , \vu \Big| r, \Theta, \vc{U} \right) =
\intO{ \left[ \frac{1}{2} \vr | \vu - \vc{U} |^2 + H_\Theta (\vr, \vt) - \frac{\partial H_\Theta (r, \Theta) }{\partial \vr} (\vr - r) -
H_\Theta (r, \Theta) \right] },
\eF
where $\{ \vr, \vt, \vu \}$ is a weak solution of the Navier-Stokes-Fourier system and $\{ r, \Theta, \vu \}$ is an arbitrary trio of smooth functions satisfying
\bFormula{i18}
r > 0, \ \Theta > 0, \ \vc{U}|_{\partial \Omega} = 0.
\eF

As observed in \cite{FeiNov10}, the weak solutions of the Navier-Stokes-Fourier system satisfy

\centerline{\textsc{Relative entropy inequality:}}
\bFormula{i19}
\left[ \mathcal{E} \left( \vr, \vt, \vu \Big| r, \Theta, \vc{U} \right) \right]_{t=0}^{t = \tau}
+ \int_0^\tau \intO{ \frac{\Theta}{\vt} \left( \tn{S}(\vt, \Grad \vu) : \Grad
\vu - \frac{ \vc{q}(\vt, \Grad \vt) \cdot \Grad \vt }{\vt} \right)} \ \dt
\eF
\[
\leq \int_0^\tau \intO{ \Big( \vr ( \vc{U} - \vu) \cdot \partial_t \vc{U} + \vr
(\vc{U} - \vu) \otimes \vu : \Grad \vc{U} - p(\vr, \vt) \Div \vc{U} \Big) } \ \dt
\]
\[
+ \int_0^\tau \intO{ \Big( \tn{S} (\vt, \Grad \vu) :
\Grad \vc{U}  } \ \dt
\]
\[
- \int_0^\tau \intO{ \left( \vr \Big( s(\vr, \vt) - s(r, \Theta) \Big) \partial_t \Theta + \vr
\Big( s(\vr, \vt) - s(r, \Theta) \Big) \vu \cdot \Grad \Theta \right) } \ \dt
\]
\[
- \int_0^\tau \intO{ \frac{\vc{q}(\vr, \Grad \vt) }{\vt} \cdot \Grad
\Theta  } \ \dt
\]
\[
+ \int_0^\tau \intO{ \left( \left( 1 - \frac{\vr}{r} \right) \partial_t p(r, \Theta) - \frac{\vr}{r} \vu \cdot
\Grad p(r, \Theta)  \right) }
\ \dt
\]
for any trio $\{ r , \Theta , \vc{U} \}$ satisfying (\ref{i18}). Taking $\{ r , \Theta, \vc{U} \}$ the strong solution of the same system emanating from the same initial data gives rise to a Gronwall type inequality for $\mathcal{E}$ yielding $\vr = r$, $\vt = \Theta$, $\vc{u} = \vc{U}$, see
\cite{FeiNov10}.

As the proof of the weak-strong uniqueness principle uses only the integral inequality (\ref{i19}) without any reference to the original system of equations, we may introduce a new class of \emph{dissipative solutions} satisfying \emph{solely} the relative entropy inequality (\ref{i19}) for any trio of
admissible \emph{test functions} $\{ r , \Theta, \vc{U} \}$. Note that a similar concept was introduced by DiPerna and Lions \cite{LI4a} in the context of the Euler system.

\subsection{Regularity criteria}

A conditional regularity criterion is a condition which, if satisfied by a weak solution, implies that the latter is regular.
Similarly, such a condition may be applied to guarantee that a local (strong) solution can be extended to a given time interval.
The most celebrated conditional regularity criteria are due to Prodi and Serrin \cite{PR}, \cite{Serrin} and Beal, Kato and Majda \cite{BeKaMa},
Constantin and Fefferman \cite{ConFef}  for solutions to the incompressible Navier-Stokes and Euler systems. Recently, similar conditions were obtained also in the context of compressible barotropic fluids
and the full Navier-Stokes-Fourier system, see
Fan, Jiang and Ou \cite{FaJiOu},  Sun et al. \cite{SuWaZh}, and the references cited therein.

In view of the results of Hoff et al. \cite{HOF7}, \cite{HofSan}, certain discontinuities imposed through the initial data in the
compressible Navier-Stokes system propagate in time. In other words, unlike its incompressible counterpart, the hyperbolic-parabolic compressible Navier-Stokes system does not enjoy the smoothing property typical for purely parabolic equations. Analogously, a solution of the full Navier-Stokes-Fourier system can be regular only if regularity is enforced by a proper choice of the initial data.

Our approach to conditional regularity is based on the concept of weak (dissipative) solutions satisfying the relative entropy inequality
(\ref{i19}):

\begin{itemize}

\item
In Section \ref{m}, we introduce the structural restrictions imposed on the thermodynamic functions $p$, $e$, $s$ and the transport
coefficients $\mu$, $\kappa$ so that the Navier-Stokes-Fourier system may possess a global in time weak (dissipative) solution
for any finite energy initial data.

\item
Now, suppose that the initial data are more regular, specifically,
\bFormula{i20}
0 < \underline{\vr} \leq \vr_0,\ 0 < \underline{\vt} \leq \vt_0 ,\ \vr_0, \vt_0 \in W^{3,2}(\Omega),
\ \vu_0 \in W^{3,2}(\Omega; R^3),
\eF
satisfying, in addition, the necessary compatibility conditions
\bFormula{i21}
\vu_0|_{\partial \Omega} = 0,\ \Grad \vt_0 \cdot \vc{n}|_{\partial \Omega} = 0, \Grad p(\vr_0, \vt_0)|_{\partial \Omega}
= \Div \tn{S}(\vt_0, \vu_0)|_{\partial \Omega}.
\eF
Under these circumstance, the Navier-Stokes-Fouries system admits a local-in-time strong solution $\{ \tilde \vr, \tilde \vt, \tilde \vu \}$ constructed by Valli \cite{Vall1} such that
\bFormula{i22}
\tilde \vr,\ \tilde \vt \in C([0,\tilde T]; W^{3,2}(\Omega)),\ \tilde \vt \in L^2(0,\tilde T; W^{4,2}(\Omega)),
\partial_t \tilde \vt \in L^2(0, \tilde T; W^{2,2}(\Omega)),
\eF
\bFormula{i23}
\tilde \vu  \in C([0,\tilde T]; W^{3,2}(\Omega;R^3)) \cap L^2(0, \tilde T; W^{4,2}(\Omega;R^3)),\
\partial_t \tilde \vu \in L^2(0,\tilde T; W^{2,2}(\Omega;R^3)).
\eF

\item Modifying slightly the arguments of \cite{FeiNov10} to accommodate the strong solutions belonging to the regularity class (\ref{i22})
we show in Section \ref{ws} that the associated (global-in-time) weak solution $\{ \vr, \vt, \vu \}$ emanating from the initial data coincides with the strong solution $\{ \tilde \vr, \tilde \vt, \tilde \vu \}$ on its existence interval $[0, \tilde T]$.

\item We show that if
\bFormula{i24}
\limsup_{t \to \tilde T-} \left( \| \vr(t, \cdot) \|_{W^{3,2}(\Omega)} + \| \vt(t, \cdot) \|_{W^{3,2}(\Omega)} + \| \vu(t, \cdot) \|_{W^{3,2}(\Omega;R^3)} \right) = \infty,
\eF
then, necessarily,
\bFormula{i25}
\limsup_{t \to \tilde T-} \| \Grad \vu (t, \cdot) \|_{L^\infty(\Omega;R^{3 \times 3})} = \infty,
\eF
see Section \ref{cr}.

\item Finally, in Section \ref{c}, we observe that the solutions constructed by Valli are classical (all necessary derivatives are continuous)
for any $t > 0$. We conclude by the following conditional regularity result for the weak solutions of the Navier-Stokes-Fourier system:

\Cbox{Cgrey}{

Suppose that $\{ \vr , \vt , \vu \}$ is a weak (dissipative) solution of the Navier-Stokes-Fourier system in $(0,T) \times \Omega$, emanating from the initial data belonging to the regularity class specified through (\ref{i20}), (\ref{i21}), and satisfying
\[
{\rm ess} \sup_{t \in (0,T)} \| \Grad \vu (t, \cdot) \|_{L^\infty(\Omega; R^{3 \times 3})} < \infty.
\]
Then $u$ is a classical solution of the Navier-Stokes-Fourier system, unique in the class of weak (dissipative) solutions,

}

\noindent see Theorem \ref{Tm1} below.
\end{itemize}

\section{Preliminaries, main results}
\label{m}

We start with a list of structural restrictions imposed on the functions $p$, $e$, $s$, $\mu$, and $\kappa$. The interested reader may consult \cite[Chapter 1]{FeNo6} for the physical background and possible relaxations. We suppose that the pressure is given by the formula
\bFormula{m1}
p(\vr, \vt) = \vt^{5/2} P \left( \frac{ \vr }{\vt^{3/2}} \right) + \frac{a}{3} \vt^4,\ a > 0,
\eF
with $P \in C^1[0,\infty) \cap C^3(0,\infty)$ satisfying
\bFormula{m2}
P(0) = 0, \ P'(Z) > 0 ,\ 0 < \frac{ \frac{5}{3} P(Z) - P'(Z) Z }{Z} < c \ \mbox{for all}\ Z > 0,
\eF
\bFormula{m3}
\lim_{Z \to \infty} \frac{ P(Z) }{Z^{5/3}}  > 0.
\eF

Accordingly, in agreement with Gibbs' equation (\ref{i11}),
\bFormula{m4}
e(\vr, \vt) = \frac{3}{2} \frac{\vt^{5/2}}{\vr} P \left( \frac{ \vr }{\vt^{3/2}} \right) + \frac{a}{\vr} \vt^4,
\eF
and
\bFormula{m5}
s(\vr, \vt) = S \left( \frac{\vr} {\vt^{3/2}} \right) + \frac{4a}{3} \frac{\vt^3}{\vr},
\eF
where
\bFormula{m6}
S'(Z) = - \frac{3}{2} \frac{ \frac{5}{3} P(Z) - P'(Z) Z }{Z^2},\ \lim_{Z \to \infty} S(Z) = 0.
\eF

Finally, the transport coefficients are continuously differentiable functions of the temperature satisfying
\bFormula{m7}
\underline{\mu} (1 + \vt^\Lambda) \leq \mu(\vt) \leq \Ov{\mu} (1 + \vt^\Lambda),\ | \mu'(\vt) | < c \ \mbox{for all}\ \vt \in [0, \infty),
\ \frac{2}{5} < \Lambda \leq 1,
\eF
\bFormula{m9}
\underline{\kappa} (1 + \vt^3) \leq \kappa (\vt) \leq \Ov{\kappa}(1 + \vt^3) \ \mbox{for all}\ \vt \in [0, \infty).
\eF

\subsection{Weak and dissipative solutions}

It was shown in \cite[Theorem 3.1]{FeNo6} that under the hypotheses (\ref{m1} - \ref{m6}), the problem
\{ (\ref{i1} - \ref{i3}), (\ref{i10}), (\ref{i13} - \ref{i15}) \} admits a global-in-time weak solution for any initial data satisfying
\bFormula{m10}
\vr_0, \vt_0 \in L^\infty(\Omega),\ \vr_0 , \vt_0 > 0 \ \mbox{a.a. in}\ \Omega,\ \vu_0 \in L^2(\Omega;R^3).
\eF
Moreover, the weak solution $\{ \vr, \vt , \vu \}$ satisfies the relative entropy inequality (\ref{i19}), see \cite{FeiNov10}.

\Cbox{Cgrey}{

Accordingly, a trio $\{ \vr, \vt, \vu \}$ will be called \emph{dissipative solution} of the Navier-Stokes-Fourier system (\ref{i1} - \ref{i8})
provided it obeys the relative entropy (\ref{i19}) for any choice of smooth functions $\{ r , \Theta , \vc{U} \}$ satisfying (\ref{i18}).

}

Summarizing the results of \cite[Chapter 3]{FeNo6} and \cite{FeiNov10} we obtain:

\bProposition{m1}
Let $\Omega \subset R^3$ be a bounded domain of class $C^{2 + \nu}$. Suppose that the thermodynamic functions $p$, $e$, $s$ and the transport coefficients $\mu$, $\kappa$ obey the structural hypotheses (\ref{m1} - \ref{m9}). Finally, let the initial data belong to the class specified in (\ref{m10}).

Then the Navier-Stokes-Fourier system possesses a dissipative solution $\{ \vr, \vt, \vu \}$ on an arbitrary time interval $(0,T)$ that enjoys the following
regularity properties:
\bFormula{m11}
\vr \geq 0 \ \mbox{a.a. in}\ (0,T) \times \Omega,\ \vr \in C([0,T]; L^1(\Omega)) \cap L^\infty(0,T; L^{5/3}(\Omega)) \cap L^\beta ((0,T) \times \Omega)
\eF
for a certain $\beta > \frac{5}{3}$;
\bFormula{m12}
\vt > 0 \ \mbox{a.a. in}\ (0,T) \times \Omega,\ \vt \in L^\infty(0,T; L^4(\Omega)) \cap L^2(0,T; W^{1,2}(\Omega)),
\eF
\bFormula{m13}
\vt^3, \log(\vt) \in L^2(0,T; W^{1,2}(\Omega));
\eF
\bFormula{m14}
\vu \in L^2(0,T; W^{\alpha}_0(\Omega;R^3)),\ \alpha = \frac{8}{5 - \Lambda},\ \vr \vu \in C_{\rm weak}(0,T; L^{5/4}(\Omega;R^3)).
\eF
\eP

\subsection{Main result}

In accordance with the programme delineated in the introductory part, our main goal is to show that for regular data the weak solution remains regular as long as we can control the gradient of the velocity field. More specifically, we will show the following result.

\Cbox{Cgrey}{

\bTheorem{m1}
Under the hypotheses of Proposition \ref{Pm1}, let $\{ \vr, \vt, \vu \}$ be a dissipative solution of the Navier-Stokes-Fourier system on the
time interval $(0,T)$ belonging to the regularity class (\ref{m11} - \ref{m14}), with the (regular) initial data satisfying (\ref{i20}), together with the compatibility conditions (\ref{i21}).

Suppose, in addition, that
\bFormula{m15}
{\rm ess} \sup_{t \in (0,T)} \| \Grad \vu (t, \cdot) \|_{L^\infty(\Omega; R^{3 \times 3})} < \infty.
\eF

Then $\{ \vr , \vt , \vu \}$ is a classical solution of the Navier-Stokes-Fourier system satisfying (\ref{i1} - \ref{i8}) in $(0,T) \times \Omega$.

\eT

}

{\bf Remark:} {\it For technical reasons, we have omitted the effect of the \emph{bulk} viscosity in the viscous stress $\tn{S} (\vt, \Grad \vu)$.
The possibility to extend the result to more general forms of the viscous stress is discussed in Section \ref{c}.}

\medskip

Here, \emph{classical solution} means that all functions and all derivatives appearing in the equations (\ref{i2} - \ref{i4}) are continuous in
$(0,T) \times \Omega$, the functions $\vr$, $\vt$, $\vu$, together with their first order derivatives, are continuous in $[0,T] \times \Omega$
and satisfy the initial conditions (\ref{i1}) as well as the boundary conditions (\ref{i7}), (\ref{i8}).

Although the rest of the paper is essentially devoted to the proof of Theorem \ref{Tm1}, we will obtain a few other results that may be of independent interest.

\section{Local existence and weak-strong uniqueness}
\label{ws}

To begin, it is convenient to introduce a scaled temperature
\[
\Xi = K(\vt), \ \mbox{where}\ K(\vt) = \int_0^\vt \kappa (z) \ {\rm d}z.
\]

Accordingly, the thermal energy balance (\ref{i12}) reads
\bFormula{ws1}
\left[ \frac{\vr}{\kappa (\vt)} \frac{\partial e(\vr, \vt)}{\partial \vt} \right] \left( \partial_t \Xi + \vu \cdot \Grad \Xi \right) - \Delta \Xi
\eF
\[
=
\left[ \frac{1}{K(\vt)}\tn{S}(\vt, \Grad \vu) : \Grad \vu \right] \Xi - \left[ \frac{\vt}{K(\vt)} \frac{\partial p(\vr, \vt)}{\partial \vt}  \right] \Div \vu\ \Xi.
\]

\subsection{First a priori bounds}
\label{WS1}

Our goal is to derive a lower and upper bounds on the temperature. We first rewrite (\ref{ws1}) in the form
\[
 \left( \partial_t \Xi + \vu \cdot \Grad \Xi \right) - D \Delta \Xi
=
D A \Xi - B \Div \vu\ \Xi,
\]
where
\[
D = \left[ \frac{\vr}{\kappa (\vt)} \frac{\partial e(\vr, \vt)}{\partial \vt} \right]^{-1},\ A =  \frac{1}{K(\vt)}\tn{S}(\vt, \Grad \vu) : \Grad \vu , \ B = \frac{\vt \kappa(\vt) }{K(\vt)} \frac{\partial p(\vr, \vt)}{\partial \vt} \left( \vr \frac{\partial e(\vr, \vt)}{\partial \vt}  \right)^{-1}.
\]

Now, in view of the hypotheses (\ref{m1}), (\ref{m4}), and (\ref{m9}), it is easy to check that the exist constants $\underline{B}$, $\Ov{B}$ such that
\[
0 < \underline{B} \leq B(t,x) \leq \Ov{B} \ \mbox{for all}\ t,x.
\]
Applying the standard comparison argument, we therefore deduce that
\bFormula{ws2}
\Xi(\tau, \cdot) \geq \inf_{x \in \Omega} \Xi(0,x) \exp \left( - \Ov{B} \int_0^\tau \| \Div \vu (t, \cdot)  \|_{L^\infty(\Omega)} \ \dt \right),\ \tau \geq 0.
\eF

In order to obtain an upper bound, we need a similar estimate for the density, namely
\bFormula{ws3}
\inf_{x \in \Omega} \vr_0(x) \exp \left( -  \int_0^\tau \| \Div \vu (t, \cdot)  \|_{L^\infty(\Omega)} \ \dt \right)
\eF
\[
\leq \vr (\tau, \cdot) \leq \sup_{x \in \Omega} \vr_0(x) \exp \left(   \int_0^\tau \| \Div \vu (t, \cdot)  \|_{L^\infty(\Omega)} \ \dt \right)
\]
that follows easily from the equation of continuity (\ref{i2}).

Now, we may use hypotheses (\ref{m7} - \ref{m9}) to observe that
\[
A (t, \cdot) \leq \Ov{A} \| \Grad \vu (t, \cdot) \|_{L^\infty(\Omega; R^{3 \times 3})}^2;
\]
whence
\bFormula{ws4}
\Xi(\tau, \cdot) \leq
\eF
\[
\sup_{x \in \Omega} \Xi(0,x)\left[ \exp \left(  \Ov{D} \
 \Ov{A} \int_0^\tau \| \Grad \vu (t, \cdot)  \|_{L^\infty(\Omega;R^{3 \times 3})}^2 \ \dt \right) + \exp \left(  \Ov{B} \int_0^\tau \| \Div \vu (t, \cdot)  \|_{L^\infty(\Omega)} \ \dt \right)
\right]
\]
for $\tau \in [0,T]$, where $\Ov{D}$ depends only on
\[
\inf_{x \in \Omega} \vr_0(x),\ \inf_{x \in \Omega} \Xi(0, x),\ \mbox{and}\ \int_0^T \| \Div \vu \|_{L^\infty(\Omega)} \ \dt .
\]

\subsection{Local existence of strong solutions}

Rewriting the system of equations (\ref{i2}), (\ref{i3}), (\ref{i12}) in terms of the new unknowns $\{ \vr, \Xi, \vu \}$
and taking the {\it a priori} bounds (\ref{ws2}), (\ref{ws4}) into account, we can apply the local existence result of
Valli \cite[Theorem A and Remark 3.3]{Vall1}. Going back to the original variables, we obtain:

\bProposition{ws1}
Under the hypotheses of Proposition \ref{Pm1}, suppose that the initial data $\{ \vr_0, \vt_0, \vu_0 \}$ satisfy
(\ref{i20}), and the compatibility conditions (\ref{i21}).

Then there exists a positive time $\tilde T$ such that the Navier-Stokes-Fourier system (\ref{i1} - \ref{i8}) admits a unique strong solution
$\{ \tvr, \tvt, \tvu \}$ in the class
\bFormula{ws5}
\tvr, \tvt \in C([0,\tilde T]; W^{3,2}(\Omega)), \ \tvu \in C([0,\tilde T]; W^{3,2}(\Omega;R^3)),
\eF
\bFormula{ws6}
\tvt \in L^2(0, \tilde T; W^{4,2}(\Omega)),\ \partial_t \tvt \in L^2(0, \tilde T; W^{2,2}(\Omega)),
\eF
\bFormula{ws7}
\tvu \in L^2(0, \tilde T; W^{4,2}(\Omega;R^3)),\ \partial_t \tvu \in L^2(0, \tilde T; W^{2,2}(\Omega;R^3)).
\eF

Moreover, there exist scalar functions $\underline{\vr}$, $\Ov{\vr}$, $\underline{\vt}$, $\Ov{\vt}$, depending solely on
\[
\inf_\Omega \vr_0, \ \sup_\Omega \vr_0, \ \inf_\Omega \vt_0, \ \sup_\Omega \vt_0, \ \mbox{and on}\
\int_0^{\tilde T} \| \Grad \tilde \vu \|^2_{L^\infty(\Omega; R^{3 \times 3})} \ \dt,
\]
such that
\bFormula{ws8}
0 < \underline{\vr}(\tau) \leq \tvr (\tau, \cdot) \leq \Ov{ \vr }(\tau),\ 0 < \underline{\vt}(\tau) \leq \tvt (\tau, \cdot) \leq \Ov{ \vt }(\tau)
\eF
for any $\tau \in [0, \tilde T]$.
\eP

\subsection{Weak strong uniqueness}

We claim that in view of the result \cite{FeiNov10} and its generalizations obtained in \cite{EF101}, the dissipative solution obtained in
Proposition \ref{Pm1} coincides with the strong solution of Proposition \ref{Pws1} on the time interval $[0, \tilde T]$ provided they start from
the same initial data (\ref{i20}), (\ref{i21}). This deserves some comments since both \cite{FeiNov10} and \cite{EF101} deal with \emph{classical solutions} having all relevant derivatives continuous and \emph{bounded} in $(0,T) \times \Omega$.

Given the integrability properties of the dissipative solutions stated in (\ref{m11} - \ref{m14}) and the regularity of the strong solutions
(\ref{ws5} - \ref{ws7}), we can easily check that the trio $\{ r = \tvr, \Theta = \tvt, \vc{U} = \tvu \}$ can be taken as test functions
in the relative entropy inequality (\ref{i19}).

Now, following step by step the arguments of \cite[Section 6, formula (79)]{EF101} we deduce from (\ref{i19}) the estimate
\bFormula{ws9}
\left[ \mathcal{E} \left( \vr, \vt, \vu \Big| \tvr, \tvt, \tvu \right) \right]_{t=0}^{t = \tau}
+ c \int_0^\tau \left[  \left\| \vu - \tilde \vu \right\|^2_{W^{1,\alpha} (\Omega;R^3)} +  \left\| \vt - \tilde \vt \right\|^2_{W^{1,2} (\Omega)}
\right] \ \dt
\eF
\[
\leq \int_0^\tau \chi_1(t) \mathcal{E} \left( \vr, \vt, \vu \Big| \tvr, \tvt, \tvu \right) \ \dt
\]
\[
+ \int_0^\tau \chi_2(t)  \intO{ \left\{ \Big( \Big[ 1 + \vr + \vr |s(\vr, \vt) | \Big]_{\rm res} \Big) \Big( \Big[ 1 +
\left| \vu - \tvu  \right| \Big]_{\rm res} \Big) + \left| \Big[  \vt - \tvt  \Big]_{\rm res} \right| \right\} }  \dt,
\]
where, similarly to \cite{EF101}, we have denoted
\[
h = [h ] _{\rm ess} + [h ]_{\rm res} = h - [h]_{\rm ess},
\]
\[
[h]_{\rm ess} = \Phi (\vr , \vt) h , \ \Phi \in \DC ((0, \infty)^2), 0 \leq \Phi \leq 1,
\]
\[
\Phi = 1
\ \mbox{in an open neighborhood of a compact}\ K \subset (0, \infty)^2
\]
where $K$ is chosen to contain the range of $[ \tvr, \tvt ]$, specifically,
\[
[\tilde \vr (t,x) , \tilde \vt (t,x) ] \in K \ \mbox{for all}\ x \in \Ov{\Omega},\ t \in [0,T].
\]
The functions $\chi_i$, $i = 1,2$ are of the form
\[
\chi_i (t)
\]
\[
= b_i(t) \Big( 1 + \| \partial_t \tvt \|_{L^\infty((0,\tilde T) \times \Omega)} + \| \nabla^2 \tvt \|_{L^\infty((0, \tilde T) \times \Omega)}
+ \| \partial_t \tvu \|_{L^\infty((0,\tilde T) \times \Omega)} + \| \nabla^2 \tvu \|_{L^\infty((0, \tilde T) \times \Omega)} \Big),
\]
where $b_i$ are bounded positive functions determined in terms of the amplitude of $\tilde \vr$, $\tilde \vt$, $\tilde \vu$ and their spatial gradients in
$[0,T] \times \Ov{\Omega}$. Focusing on the most difficult term, we have
\[
\chi_2(t)
\intO{ \vr |s(\vr,\vt)| | \left[ \vu - \tvu \right]_{\rm res} | } \leq \chi_2(t) \left\| [ \vr s(\vr, \vt) ]_{\rm res} \right\|_{L^{4/3}(\Omega)}
\left\| \vu - \tvu \right\|_{L^4(\Omega;R^3)}
\]
\[
\leq \ep \left\| \vu - \tvu \right\|^2_{W^{1,\alpha} (\Omega;R^3)} + c(\ep) \chi^2_2(t) \left\| [ \vr s(\vr, \vt) ]_{\rm res} \right\|^2_{L^{4/3}(\Omega)}
\]
\[
\leq \ep \left\| \vu - \tvu \right\|^2_{W^{1,\alpha} (\Omega;R^3)} + c(\ep) \chi_2^2(t) \mathcal{E}^{3/2} \left( \vr , \vt, \vu \Big| \tvr, \tvt, \tvu \right) \ \mbox{for any}\ \ep > 0,
\]
where, exactly as in \cite[Section 6.1]{EF101}, we have used the structural properties of $s$ and the embedding
\[
W^{1,\alpha}(\Omega) \hookrightarrow L^4,\ \alpha = \frac{8}{5 - \Lambda},\ \frac{2}{5} < \Lambda \leq 1.
\]

Going back to (\ref{ws9}) we may infer that
\bFormula{ws10}
\left[ \mathcal{E} \left( \vr, \vt, \vu \Big| \tvr, \tvt, \tvu \right) \right]_{t=0}^{t = \tau}
+ c \int_0^\tau \left[  \left\| \vu - \tilde \vu \right\|^2_{W^{1,\alpha} (\Omega;R^3)} +  \left\| \vt - \tilde \vt \right\|^2_{W^{1,2} (\Omega)}
\right] \ \dt
\eF
\[
\leq \int_0^\tau \chi(t) \mathcal{E} \left( \vr, \vt, \vu \Big| \tvr, \tvt, \tvu \right) \ \dt,
\]
where, by virtue of (\ref{ws6}), (\ref{ws7}),
\[
\chi \in L^1(0,T).
\]
Applying Gronwall's lemma we obtain the following conclusion:

\bProposition{ws2}
Under the hypotheses of Proposition \ref{Pws1}, let $\{\tvr, \tvt, \tvu \}$ be the local strong solution specified in Proposition
\ref{Pws1} and $\{ \vr, \vt , \vu \}$ a dissipative solution obtained in Proposition \ref{Pm1}, emanating from the same initial data.

Then
\[
\vr = \tvr,\ \vt = \tvt,\ \vu = \tvu \ \mbox{in}\ [0, \tilde T] \times \Omega.
\]
\eP

\section{Conditional regularity}
\label{cr}

Our goal is to show that the energy norm
\[
 \| \tvr(t, \cdot) \|_{W^{3,2}(\Omega)} + \| \tvt(t, \cdot) \|_{W^{3,2}(\Omega)} + \| \tvu(t, \cdot) \|_{W^{3,2}(\Omega;R^3)}
\]
associated to a strong solution of the Navier-Stokes-Fourier system remains bounded in $[0,\tilde T]$ as long as
\bFormula{cr1}
\sup_{t \in (0, \tilde T)} \| \Grad \tilde \vu (t, \cdot) \|_{L^\infty(\Omega;R^{3 \times 3})} \leq G < \infty.
\eF
We remark that we already know that
\bFormula{cr2}
0 < \underline{\vr}(\tau) \leq \tvr (\tau, \cdot) \leq \Ov{ \vr }(\tau),\ 0 < \underline{\vt} \leq \tvt (\tau, \cdot) \leq \Ov{ \vt }
\eF
for any $\tau \in [0, \tilde T]$, see (\ref{ws8}), where the bounds depend only on $G$, the initial data, and the length of the time interval.

\subsection{Energy bounds, temperature}

Multiplying equation (\ref{ws1}) on
\[
\left[ \frac{\tvr}{\kappa (\tvt)} \frac{\partial e(\tvr, \tvt)}{\partial \vt} \right]^{-1} \Delta \tilde \Xi,\ \tilde \Xi = K(\tvt),
\]
we obtain, after a routine manipulation,
\bFormula{cr3}
{\rm ess} \sup_{t \in (0, \tilde T)} \| \tilde \Xi (t, \cdot) \|_{W^{1,2}(\Omega)} +
\left\| \partial_t \tilde \Xi \right\|^2_{L^2(0,T; L^2(\Omega))} + \left\| \Xi \right\|^2_{L^2(0,T; W^{2,2}(\Omega))} \leq c(B, {\rm data}).
\eF

Since $\tvt$ is already known to be bounded, the estimates (\ref{cr3}) transfer to $\tvt$. Indeed note that
\[
D^2_x \tvt = \left[ K^{-1}(\tilde \Xi ) \right]^{-1} D_x^2 \tilde \Xi + \left[ K^{-1}(\tilde \Xi) \right] |D_x \tilde \Xi |^2,
\]
where, since $\tilde \Xi$ satisfies the homogeneous Neumann boundary condition, the term $\Grad \tilde \Xi$ may be estimates in terms of the
Gagliardo-Nirenberg inequality
\bFormula{cr4}
\| \Grad \tilde \Xi \|^2_{L^4(\Omega;R^3)} \leq \| \tilde \Xi \|_{L^\infty(\Omega)} \| \Delta \Xi \|_{L^2(\Omega)}.
\eF
Consequently, relation (\ref{cr3}) implies
\bFormula{cr5}
{\rm ess} \sup_{t \in (0, \tilde T)} \| \tvt (t, \cdot) \|_{W^{1,2}(\Omega)} +
\left\| \partial_t \tvt \right\|^2_{L^2(0,T; L^2(\Omega))} + \left\| \tvt \right\|^2_{L^2(0,T; W^{2,2}(\Omega))} \leq c(B, {\rm data}).
\eF

\subsection{Energy bounds, velocity}

Our next goal is to deduce similar bounds for the velocity field. To this end, we write the momentum equation in the form
\bFormula{cr6}
\partial_t \tvu - \frac{1}{\tvr} \Div \tn{S} (\tvt, \Grad \tvu) = - \tvu \cdot \Grad \tvu - \frac{1}{\tvr} \frac{\partial p(\tvr, \tvt) }{\partial \vt} \Grad \tvt + \frac{1}{\tvr} \frac{\partial p(\tvr, \tvt) }{\partial \vr} \Grad \tvr
\eF
\[
= \vc{h}_1 + \frac{1}{\tvr} \frac{\partial p(\tvr, \tvt) }{\partial \vr} \Grad \tvr,
\]
where, in accordance with the previous estimates
\[
\sup_{t \in [0, \tilde T] } \| \vc{h}_1(t,\cdot) \|_{L^2(\Omega;R^3)} \leq c(B, {\rm data}).
\]

Taking the scalar product of (\ref{cr6}) with $- \Div \tn{S} (\tvt, \Grad \tvu)$ we obtain
\[
\intO{ \tn{S} (\tvt, \tvu) : \partial_t \Grad \tvu }  + \intO{ \frac{1}{\tvr} \left| \Div \tn{S} (\tvt, \Grad \tvu) \right|^2 }
\]
\[
= - \intO{ \left( \vc{h}_1 + \frac{1}{\tvr} \frac{\partial p(\tvr, \tvt) }{\partial \vr} \Grad \tvr \right) \cdot \Div \tn{S} (\tvt, \Grad \tvu) },
\]
where, furthermore,
\[
\intO{ \tn{S} (\tvt, \tvu) : \partial_t \Grad \tvu }
\]
\[
= \frac{{\rm d}}{{\rm d}t}  \intO{  \frac{\mu(\tvt)}{4}
\left| \Grad \tvu + \Grad^t \tvu - \frac{2}{3} \Div \tvu \tn{I} \right|^2 }
\]
\[
- \intO{ \frac{\mu'(\tvt)}{4} \partial_t \tvt \left| \Grad \tvu + \Grad^t \tvu - \frac{2}{3} \Div \tvu \tn{I} \right|^2 } .
\]

In view of the estimate (\ref{cr5}) we conclude that
\bFormula{cr7}
\frac{{\rm d}}{{\rm d}t}  \intO{  \frac{\mu(\tvt)}{4}
\left| \Grad \tvu + \Grad^t \tvu - \frac{2}{3} \Div \tvu \tn{I} \right|^2 }
+ \intO{ \frac{1}{\tvr} \left| \Div \tn{S} (\tvt, \Grad \tvu) \right|^2 }
\eF
\[
= - \intO{ \left( \vc{h}_1 + \frac{1}{\tvr} \frac{\partial p(\tvr, \tvt) }{\partial \vr} \Grad \tvr \right) \cdot \Div \tn{S} (\tvt, \Grad \tvu) }
+ \intO{ h_2 },
\]
where
\[
\| \vc{h}_2  \|_{L^2(0,\tilde T ; L^2(\Omega)} \leq c(B, {\rm data}).
\]

\subsection{Elliptic estimates}\label{ee}

In order to exploit (\ref{cr7}), we have to show that the $L^2$-norm of $\Div \tn{S}(\tvt, \Grad \tvu)$ is ``equivalent'' to
the $L^2$-norm of $\nabla^2 \tvu$. We have
\[
\vc{H} \equiv \Div \tn{S}(\tvt, \Grad \tvu) = \mu(\tvt) \left[ \left( \Delta \tvu + \frac{1}{3} \Grad \Div \tvu \right) \right] + \vc{h}_3,
\]
where, in accordance with (\ref{cr5}),
\[
\| \vc{h}_3  \|_{L^\infty(0,\tilde T ; L^2(\Omega;R^3)} \leq c(B, {\rm data}).
\]
Consequently, we may integrate (\ref{cr7}) with respect to $t$ and use the standard elliptic estimates to obtain
\bFormula{cr8}
\int_0^\tau \| \tvu \|^2_{W^{2,2}(\Omega;R^3)} \ \dt \leq c(B, {\rm data}) \left( 1 + \int_0^\tau \| \Grad \tvr \|^2_{L^2(\Omega;R^3)} \right),
\ \tau \in [0, \tilde T].
\eF

Now, we differentiate the equation of continuity with respect to $x$ to obtain
\bFormula{cr9}
\partial_t \left( \partial_{x_i} \tvr \right) + \tvu \cdot \Grad \left( \partial_{x_i} \tvr \right) =
- \partial_{x_i}(\tvu) \cdot \Grad \tvr - \left( \partial_{x_i} \tvr \right) \Div \tvu - \tvr \partial_{x_i} \Div \tvu;
\eF
whence
\bFormula{cr10}
\frac{{\rm d}}{{\rm d}t} \intO{ |\Grad \tvr |^2 } \leq c(B, {\rm data}) \intO{ \left( | \Grad \tvr |^2 + |\Grad \tvr || \Grad \Div \tvu |
\right)}.
\eF

Combining (\ref{cr8}), (\ref{cr10}) with the standard Gronwall type argument we get the following estimates:
\bFormula{cr11}
\sup_{t \in [0, \tilde T]} \| \Grad \tvr (t, \cdot) \|_{L^2(\Omega;R^3)} \leq c(B, {\rm data}).
\eF
Next, due to (\ref{cr6}), 
\bFormula{cr12}
\left\| \partial_t \tvu \right\|_{L^2(0,\tilde T; L^2(\Omega;R^3))} + \left\|  \tvu \right\|_{L^2(0,\tilde T; W^{2,2}(\Omega;R^3))} \leq c(B, {\rm data}),
\eF
and, finally, by means of the equation of continuity,
\bFormula{cr13}
\left\| \partial_t \tvr \right\|_{L^2(0,\tilde T; L^2(\Omega))} \leq c(B, {\rm data}).
\eF

\subsection{Energy estimates for the time derivatives}

Differentiating the momentum equation (\ref{i3}) with respect to $t$ and setting $\tvV = \partial_t \tvu$ we obtain
\bFormula{cr14}
\tvr \left( \partial_t \tvV + \tvu \cdot \Grad \tvV \right) - \Div  \tn{S}(\tvt, \Grad \tvV)
= - \left( \partial_t \tvr \tvV + \tvr \tvV \cdot \Grad \tvu \right)  + \vc{h}_1,
\eF
with
\[
\vc{h}_1 = - \partial_t \tvr \tvu \cdot \Grad \tvu + \Div \left[ \left( \mu'(\tvt) \partial_t \tvt \right)\left( \Grad \tvu +
\Grad^t \tvu - \frac{2}{3} \Div \tvu \tn{I} \right) \right] - \Grad \left( \partial_t p(\tvr, \tvt) \right),
\]
where, in accordance with the previous estimates,
\bFormula{cr15}
\| \vc{h}_1 \|_{L^2(0,\tilde T; W^{-1,2}(\Omega;R^3))} \leq c(G, {\rm data}).
\eF

Seeing that
\[
- \left( \partial_t \tvr \tvV + \tvr \tvV \cdot \Grad \tvu \right) \cdot \vc{V} = |\tvV |^2 \tvr \Div \tvu
\]
we may take the scalar product of (\ref{cr14}) with $\tvV$ and, integrating over $\Omega$, we deduce the energy estimates for
$\tvV = \partial_t \tvu$:
\bFormula{cr16}
\sup_{t \in [0,\tilde T]} \| \partial_t \tvu \|_{L^2(\Omega;R^3)} + \| \partial_t \tvu \|_{L^2(0,\tilde T; W^{1,2}(\Omega;R^3))} \leq
c(G, {\rm data}).
\eF
and, according to Sobolev's embedding theorem,
\bFormula{cr22}
\left\| \partial_t \tvu \right\|_{L^2(0, \tilde T; L^{6}(\Omega;R^3))} \leq c(G, {\rm data}).
\eF

Going back to (\ref{i3}) we compute
\bFormula{cr17}
\mu(\tvt) \left( \Delta \tvu + \frac{1}{3} \Grad \Div \tvu \right)
\eF
\[
= \tvr \left( \partial_t \tvu + \tvu \cdot \Grad \tvu \right) -
\mu'(\tvt) \left( \Grad \tvu + \Grad \tvu^t - \frac{2}{3} \Div \tvu \tn{I} \right) \Grad \tvt + \Grad p(\tvr, \tvt);
\]
which, combined with (\ref{cr16}) and the previous estimates, gives rise to
\bFormula{cr18}
\sup_{t \in [0, \tilde T]} \| \tvu (t, \cdot) \|_{W^{2,2}(\Omega;R^3)} \leq  c(G, {\rm data}).
\eF

Next, bootstraping (\ref{cr17}) via the elliptic regularity,  yields
\bFormula{cr19}
\int_0^\tau \| \tvu(t, \cdot) \|^2_{W^{2,6}(\Omega; R^3)} \ \dt \leq c \left( 1 +  \int_0^\tau \| \Grad \tvr (t, \cdot) \|_{L^6(\Omega;R^3)}^2
\ \dt \right).
\eF
Furthermore, multiplying (\ref{cr9}) on $| \Grad \tvr |^4 \Grad \tvr$ yields
\bFormula{cr20}
\frac{{\rm d}}{{\rm d}t} \| \Grad \tvr \|^6_{L^6(\Omega;R^3)} \leq c \left( \| \Grad \tvr \|^6_{L^6(\Omega;R^3)} + \| \tvu \|_{W^{2,6}(\Omega;R^3)}
\| \Grad \tvr \|_{L^6(\Omega;R^3)}^{5} \right).
\eF
Thus, combining (\ref{cr19}), (\ref{cr20}), we may infer that
\bFormula{cr21}
\sup_{t \in [0, \tilde T]} \left[ \| \Grad \tvr (t, \cdot) \|_{L^6(\Omega;R^3)} + \| \partial_t \tvr (t, \cdot) \|_{L^6(\Omega;R^3)} \right]
+ \left\| \tvu \right\|_{L^2(0, \tilde T; W^{2,6}(\Omega;R^3))} \leq c(G, {\rm data}).
\eF

\subsection{$L^p - L^q$ estimates for the temperature}

Our next goal is to apply the technique of $L^p - L^q$ estimates to the parabolic equation (\ref{i12}). For this purpose, we first need H\"older continuity for $\tvt$. To this end, let us write equation (\ref{i12}) in terms of $\tilde e(\tvr, \tilde\Xi) = e(\tvr, K^{-1}(\tilde \Xi))$.
\[
\tvr\partial_t \tilde e(\tvr,\tilde\Xi) + \tvr \tvu \cdot \Grad \tilde e(\tvr,\tilde\Xi) -  \Delta_x \tilde\Xi = \tn{S}(\tvt, \Grad \tvu) : \Grad \tvu  - p(\tvr, \tvt) \Div \tvu := h.
\]
Note that
\[
\Delta_x \tilde\Xi  = \Div \left( \frac{\Grad e }{e_{\Theta}(\tvr, \tilde\Xi)} \right) - \Div \left( \frac{ \partial_\vr \tilde e (\tvr, \tilde\Xi) }{\partial_\Theta \tilde e(\tvr, \tilde\Xi)} \Grad \tvr \right).
\]
We obtain
\[
\tvr\partial_t \tilde e(\tvr,\tilde\Xi) + \tvr \tvu \cdot \Grad \tilde e(\tvr,\tilde\Xi) -  \Div \left( \frac{\Grad 
\tilde e }{\partial_\Theta \tilde e (\tvr, \tilde\Xi)} \right) = h + \Div \left( \frac{ \partial_{\vr} \tilde e_(\tvr, \tilde\Xi) }{ \partial_\Theta \tilde e (\tvr, \tilde\Xi)} \Grad \tvr \right),
\]
which, after dividing on both sides by $\tvr$, yields
\[
\partial_t \tilde e(\tvr,\tilde\Xi) + \left( \tvu - \frac{\Grad \tvr}{\tvr^2 { \partial_\Theta \tilde e (\tvr, \tilde\Xi)}} \right) \cdot \Grad \tilde e(\tvr,\tilde\Xi)  -  \Div \left( \frac{\Grad e }{\tvr \partial_\Theta \tilde e (\tvr, \tilde\Xi)} \right)
\]
\bFormula{cr23}
= \frac{h}{\tvr}+ \frac{ \partial_\vr \tilde e (\tvr, \tilde\Xi) }{ \tvr^2 \partial_\Theta \tilde e (\tvr, \tilde\Xi)}|\Grad \tvr|^2 + \Div \left( \frac{ \partial_\vr \tilde e (\tvr, \tilde\Xi) }{ \tvr \partial_\Theta \tilde e (\tvr, \tilde\Xi)} \Grad \tvr \right).
\eF
Note that $\Grad \tvr \in L^{\infty}(0,\tilde T; L^6 (\Omega, R^3))$ according to (\ref{cr21}), we can apply the standard theory of parabolic equations with bounded measurable coefficients to deduce that
\[
\tilde e(\tvr, \tilde\Xi) \ \mbox{is H\" older continuous in}\ [0,T] \times \Ov{\Omega}.
\]
See Ladyzhenskaya et al. \cite{LADSOLUR}. Since $\tvr$ is already H\" older continuous
in the set $[0,T] \times \Ov{\Omega}$ (cf. the estimates (\ref{cr21})), we find
\bFormula{cr24}
\tilde \Xi (\mbox{hence}\, \tvt) \ \mbox{is H\" older continuous in}\ [0,T] \times \Ov{\Omega}.
\eF
Now we write (\ref{i12}) in the following form:
\bFormula{cr25}
\partial_t \tilde \Xi + \tvu\cdot \Grad \tilde\Xi - D \Delta_x \tilde \Xi = D h.
\eF
with the diffusion coefficient $D = \left( \tvr \frac{\partial e(\tvr, \tvt)}{\partial \vt} \right)^{-1}$ being H\" older continuous. The $L^p-L^q$ theory for parabolic equations is now applicable yielding
\bFormula{cr26}
\tvt \in L^p(0,\tilde T; W^{2,6}(\Omega)),\ \partial_t \tvt \in L^p(0,\tilde T; L^6(\Omega)) \ \mbox{for any}\ 1 < p < \infty.
\eF
see Amann  \cite{Amann1}, \cite{AMA}, Krylov \cite{Krylov}. Note that the integrability in the spatial variable is limited by the
integrability of the initial data.

{\bf Remark: }
Applying similar arguments to the momentum equation (\ref{i3}), we could obtain analogous estimates for the velocity field:
\bFormula{cr27}
\tvu \in L^p(0,\tilde{T}; W^{2,6}(\Omega;R^3)),\ \partial_t \tvu \in L^p(0,\tilde T; L^6(\Omega;R^3)) \ \mbox{for any}\ 1 < p < \infty,
\eF
however, we do not need this refinement in the future analysis.

\subsection{Full regularity}

Differentiating (\ref{cr25}) with respect to time yields
\bFormula{cr29}
\partial_t \tilde \Xi_t  - D \Delta_x \tilde \Xi_t = (D h)_t + D_t \Delta_x \tilde \Xi  - (\tvu\cdot \Grad \tilde\Xi)_t := \tilde h .
\eF
According to (\ref{cr16}), (\ref{cr21}) and (\ref{cr27}), it is easy to check that
\bFormula{cr30}
\| \tilde{h}\|_{L^2(0, \tilde T; L^2(\Omega))} \leq c(G, {\rm data}).
\eF
Standard energy method and elliptic estimates yield
\bFormula{cr30a}
\| \tilde \Xi_t \|_{L^{\infty}(0,\tilde T; L^2(\Omega))} + \| \tilde \Xi_t \|_{L^{2}(0, \tilde T; W^{1,2}(\Omega))} \leq c(G, {\rm data}),
\eF
\bFormula{cr31}
\| \tilde \Xi \|_{L^{\infty}(0,\tilde T; W^{2,2}(\Omega))} + \| \tilde \Xi \|_{L^{2}(0,\tilde T; W^{3,2}(\Omega))} \leq c(G, {\rm data}),
\eF
as well as
\bFormula{cr32}
\| \tilde \Xi_t \|_{L^{\infty}(0, \tilde T; W^{1,2}(\Omega))} + \| \tilde \Xi_{tt} \|_{L^{2}(0,\tilde T; L^2(\Omega))} \leq c(G, {\rm data}),
\eF
\bFormula{cr33}
\| \tilde \Xi \|_{L^{\infty}(0,\tilde T; W^{3,2}(\Omega))} + \| \tilde \Xi_t \|_{L^{2}(0,\tilde T; W^{2,2}(\Omega))} \leq c(G, {\rm data}),
\eF
\bFormula{cr34}
\| \tilde \Xi \|_{L^{2}(0,\tilde T; W^{4,2}(\Omega))} \leq c(G, {\rm data}).
\eF
Similar estimates for $\tvt$ hold. With these estimates in hand, we can go back to equation (\ref{cr14}) for $\tilde\vc{V} = \partial_t \tvu$ and
in a similar way to conclude
\bFormula{cr35}
\| \tvu_t \|_{L^{\infty}(0, \tilde T; W^{1,2}(\Omega,R^3))} + \| \tvu_{tt} \|_{L^{2}(0,\tilde T; L^2(\Omega,R^3))} \leq c(G, {\rm data}),
\eF
\bFormula{cr36}
\| \tvu \|_{L^{\infty}(0,\tilde T; W^{3,2}(\Omega,R^3))} + \| \tvu_t \|_{L^{2}(0,\tilde T; W^{2,2}(\Omega,R^3))} \leq c(G, {\rm data}).
\eF
As for the final estimates
\bFormula{cr37}
\| \tvu \|_{L^{2}(0,\tilde T; W^{4,2}(\Omega,R^3))} + \| \tvr \|_{L^{\infty}(0, \tilde T; W^{3,2}(\Omega))} \leq c(G, {\rm data}),
\eF
one needs to combine with the transport equation. The treatment is similar to Section \ref{ee} and we omit the details.

Summarizing the previous considerations, we are allowed to state the following result.

\bProposition{cr1}
Let $T > 0$ be given.
Suppose that $\{ \tvr, \tvt, \tvu \}$ is a strong solution of the Navier-Stokes-Fourier on the time interval $[0, \tilde T]$, $
\tilde T \leq T$,  the existence
of which is claimed in Proposition \ref{Pws1}. Assume, in addition, that
\bFormula{cr38}
\sup_{t \in [0, \tilde T]} \| \Grad \tvu (t, \cdot) \|_{L^\infty(\Omega;R^3)} \leq G.
\eF

Then
\bFormula{cr39}
\sup_{t \in [0, \tilde T]} \left[ \| \tvr(t,\cdot) \|_{W^{3,2}(\Omega)} + \| \tvt(t,\cdot) \|_{W^{3,2}(\Omega)} +
\| \tvu (t, \cdot) \|_{W^{3,2}(\Omega;R^3)} \right] \leq c(G, {\rm data}, T).
\eF
In particular, if $\tilde T$ is the \emph{maximal} existence time, then $\tilde T = T$.
\eP

\section{Conclusion}

\label{c}

To begin, the standard parabolic regularity theory implies that the solutions belonging to the class specified in Proposition \ref{Pws1}
are, in fact, classical solutions, meaning all relevant derivatives appearing in the system (\ref{i2} - \ref{i4}) are continuous in the
\emph{open} set $(0,\tilde T) \times \Omega$, the fields $\tvt$, $\tvr$ are continuously differentiable up to the boundary $\partial \Omega$, and
the boundary conditions (\ref{i7}), (\ref{i8}) hold, see Matsumura and Nishida \cite{MANI1}, \cite{MANI}. Consequently, combining Proposition \ref{Pws1} with Proposition \ref{Pcr1} we obtain the conclusion of Theorem \ref{Tm1}.

We conclude the paper by comments concerning the hypotheses of Theorem \ref{Tm1}. In the light of the results obtained by Fan et al. \cite{FaJiOu},
the ``optimal'' regularity criterion is expected to be the bound
\[
\int_0^T \| \Grad \vu \|_{L^\infty(\Omega;R^{3 \times 3})} \ \dt < \infty.
\]
Note, however, that, unlike \cite{FaJiOu}, we have to handle the transport coefficients that depend effectively on the temperature. The advantage of this approach is that, again unlike \cite{FaJiOu}, the Navier-Stokes-Fourier system is known to possess a global in time dissipative solution.

The regularity of the strong solutions seems optimal, at least in view of the result \cite{FeiNov10}. Indeed the time derivative $\partial_t \tvu$ of
the strong solutions should be at least of class $L^2(0,T; L^\infty(\Omega;R^3))$ for the technique of \cite{FeiNov10} to be applicable.

The next remark concerns the possibility to include the more general class of Newtonian stress, namely
\[
\tn{S}(\vt, \Grad \vu) = \mu(\vt) \left( \Grad \vu + \Grad^t \vu - \frac{2}{3} \Div \vu \tn{I} \right) + \eta (\vt) \Div \vu \tn{I},
\]
where we have added the bulk viscosity, with the coefficient $\eta = \eta(\vt)$. Clearly, we can show the same result if
\[
\frac{\eta(\vt)}{\mu(\vt)} = A  - \mbox{a constant.}
\]
Moreover, our method still applies if
\[
\frac{\eta(\vt)}{\mu(\vt)} \approx \mbox{``small''},
\]
where \emph{small} is determined in terms of the optimal constant $\lambda$ in the elliptic estimate
\[
\left\| \Delta \vc{v} + \frac{1}{3} \Div \vc{v} \tn{I} \right\|_{L^2(\Omega; R^{ 3 \times 3})} \geq \lambda \| \vc{v} \|_{W^{2,2}(\Omega;R^3)},\
\vc{v}|_{\partial \Omega} = 0.
\]

The general case may be handled in the following way:

\begin{enumerate}

\item After establishing the uniform bounds in Section \ref{WS1}, apply the theory of Krylov and Safonov \cite{KrySa} for parabolic equation in 
\emph{non-divergence} form to equation (\ref{ws1}) to obtain uniform bounds on $\tvt$ in the H\" olderian norm.

\item Use the regularity theory for elliptic systems with H\" older coefficients in Section \ref{ee}.

\end{enumerate}

Finally, we remark that the methods of the present paper could be extended to other types of boundary conditions of the Navier type for the velocity and to more general classes of domains $\Omega$ including certain unbounded domains like exterior domains or a half-space.


\end{document}